\documentclass[12pt]{article}
\usepackage{amsmath,amssymb,latexsym,epsfig}
\newcommand{\beq}[1]{\begin{equation}\label{#1}}
\newcommand{\eeq}{\end{equation}}

\newcommand{\set}[1]{\left\{#1\right\}}

\newcommand{\ignore}[1]{}

\def\cE{{\cal E}}

\def\cP{{\cal P}}

\def\e{\epsilon}

\def\G{\Gamma}

\def\bX{\bar{X}}

\def\Pr{\mbox{{\bf Pr}}}

\def\whp{{\bf whp}}

\newtheorem{theorem}{Theorem}

\newcommand{\Hp}[1]{H_{n,#1;3}}
\def\cY{{\cal Y}}

\def\cP{{\cal P}}

\begin{document}
\title{{\bf Loose Hamilton Cycles in Random 3-Uniform Hypergraphs}}
\author{Alan Frieze\thanks{
Department of Mathematical Sciences, Carnegie Mellon University,
Pittsburgh PA15213, U.S.A. Supported in part by NSF grant
DMS-0753472.} }
\maketitle
\begin{abstract}
In the random hypergraph $H=H_{n,p;3}$ each possible triple appears independently 
with probability $p$. A {\em loose} Hamilton cycle can be described as a sequence of edges
$\set{x_i,y_i,x_{i+1}}$ for $i=1,2,\ldots,n/2$. We prove that
there exists an absolute constant $K>0$ such that if $p\geq \frac{K\log n}{n^2}$ then
$$\lim_{\substack{n\to \infty\\4 |n}}\Pr(H_{n,p;3}\ contains\ a\ loose\ Hamilton\ cycle)=1.$$
\end{abstract}

\section{Introduction}
The threshold for the existence of Hamilton cycles in the random graph $G_{n,p}$ has been known
for many years, see \cite{KS}, \cite{AKS} and \cite{Boll}. There have been many generalisations
of these results over the years and the problem is well understood. It is natural to try to extend
these results to Hypergraphs and this has proven to be difficult. The famous P\'osa lemma fails to 
provide any comfort and we must seek new tools. In the graphical case, Hamilton cycles and perfect
matchings go together and our approach will be to build on the deep and difficult result of
Johansson, Kahn and Vu \cite{JKV}, as well as what we have learned from the graphical case.

A $k$-uniform Hypergraph is a pair $H=(V,E)$ where $E\subseteq \binom{V}{k}$.  We say that a $k$-uniform 
sub-hypergraph $C$ of $H$ is a Hamilton cycle of type
$\ell$, for some $1\le \ell \le k$, if there exists a cyclic
ordering of the vertices $V$ such that every edge consists of $k$
consecutive vertices and for every pair of consecutive edges
$E_{i-1},E_i$ in $C$ (in the natural ordering of the edges) we have
$|E_{i-1}\setminus E_i|=\ell$. When $\ell=k-1$ we say that $C$ is a {\em loose} Hamilton cycle and in this
paper we will restrict our attention to loose Hamilton cycles in the random 3-uniform 
hypergraph $H=H_{n,p;3}$. In this hypergraph, $V=[n]$ and 
each of the $\binom{n}{3}$ possible edges (triples) appears independently with probability $p$.
While $n$ needs to be even for $H$ to contain a loose Hamilton cycle, we need to go one step further
and assume that $n$ is a multiple of 4.
Extensions to other $k,\ell$ and $n=2\mod 4$ pose problems. We will prove the following theorem:
\begin{theorem}\label{th1}
There exists an absolute constant $K>0$ such that if $p\geq \frac{K\log n}{n^2}$ then
$$\lim_{\substack{n\to \infty\\4 |n}}\Pr(H_{n,p;3}\ contains\ a\ loose\ Hamilton\ cycle)=1.$$
\end{theorem}
Thus $\frac{\log n}{n^2}$ is the threshold for the existence of loose Hamilton cycles, at least for 
$n$ a multiple of 4. This is because
if $p\leq \frac{(1-\e)\log n}{2n^2}$ and $\e>0$ is constant, then \whp\footnote{An event $\cE_n$
occurs {\em with high probability}, or \whp\ for brevity, if
$\lim_{n\rightarrow\infty}\Pr(\cE_n)=1$.} $\Hp{p}$ contains isolated vertices.

The proof of Theorem \ref{th1} will follow fairly easily from the following three theorems. 

We start with a special case of the theorem of \cite{JKV}: 
Let $X$ and $Y$ be a disjoint sets.
Let $\Omega=\binom{X}{2}\times Y$.
Let $\G=\G(X,Y,p)$ be the random 3-uniform hypergraph where each triple in $\Omega$ is independently
included with probability $p$. 
Assuming that $|X|=2|Y|=2m$, a perfect matching of $\G$ is a set of $m$ triples $(x_{2i-1},y_i,x_{2i}),\,i=1,2,\ldots,m$
such that $X=\set{x_1,\ldots,x_{2m}}$ and $Y=\set{y_1,\ldots,y_m}$. 
\begin{theorem}\label{th2}\cite{JKV}\ \\
There exists an absolute constant $K>0$ such that if $p\geq \frac{K\log n}{n^2}$ then 
\whp\ $\G$ contains a perfect matching.
\end{theorem}
This version is not actually proved in \cite{JKV}, but can be obtained by straightforward changes to their proof.

Our next theorem concerns {\em rainbow} Hamilton cycles in random regular graphs. If we edge colour a graph then
a set $S$ of edges is rainbow if all edges in $S$ are a different colour. Janson and Wormald \cite{JW} proved the following.
\begin{theorem}\label{th3}
If the edges of a random $2r$-regular graph $G_{2r}$ on vertex set $[n]$ are coloured randomly with
$n$ colours so that each colour is used exactly $r$ times, $r\geq 4$, then \whp\ it contains a rainbow Hamilton
cycle.  
\end{theorem}

The relationship between a loose Hamilton cycle $(x_1,y_1,x_2),(x_2,y_2,x_3),\ldots,(x_m,y_m,x_1)$ becomes apparent
if we consider $y_1,y_2,\ldots,y_m$ to be the colors of edges $(x_i,x_{i+1}),i=1,2,\dots,m$ in an associated graph.
More precisely, we will partition $[n=4m]$ into $X=[2m]$ and $\bX=[2m+1,n]$. The 
(multi-)graph $G^*$ has vertex set $X$ and an edge 
$(x,x')$ of colour $y$ if $(x,y,x')$ is an edge of $H$. If $G^*$ contains a rainbow Hamilton cycle, then $H$
contains a loose Hamilton cycle. We will use Theorem \ref{th2} to show that \whp\ $G^*$ contains an edge coloured
graph that is close to satisfying the conditions of Theorem \ref{th3}.

There is a minor technical point in that we can only use Theorem \ref{th2} to prove the existence
of a randomly coloured (multi-)graph $\G_{2r}$ that is the union of $2r$ independent matchings. Fortunately, 
\begin{theorem}\label{th4}
$\G_{2r}$ is {\em contiguous} to $G_{2r}$  
\end{theorem}
By this we mean that if $\cP_n$ is some sequence of (multi-)graph properties, then
\beq{cP}
\G_{2r}\in\cP_n\ \whp\ \Longleftrightarrow\ G_{2r}\in\cP_n\ \whp.
\eeq
Theorem \ref{th4} is proved in Janson \cite{J} and in Molloy, Robalewska-Szalat, Robinson and Wormald \cite{MRRW}.
\section{Proof of Theorem \ref{th1}}
We begin by letting $\cY$ be a set of size $2rm$ consisting of $r=O(1)$ copies $y_1,y_2,\ldots,y_r$ of each $y\in \bX$.
Next let $Y_1,Y_2,\ldots,Y_{2r}$ be a uniformly random partition of $\cY$ into $2r$ sets of size $m$.

Define $p_1$ by $p=1-(1-p_1)^{2r}$. With this choice, we can generate $H_{n,p;3}$ as the union of $2r$ independent copies of 
$H_{n,p_1;3}$. Similarly, define $p_2$ by $p_1=1-(1-p_2)^r$.

Viewing $H_{n,p_1;3}$ as the union of $r$ independent copies $H_1,H_2,\ldots,H_r$ of 
$H_{n,p_2;3}$ we can couple $\G(X,Y_j,p_1)$ with a subgraph of $H_{n,p_1;3}$ by placing $(x,y,x')$ in $E(H_i)$ if 
$(x,y_i,x')\in E(\G(X,Y_i,p_1))$. 

It follows from Theorem \ref{th2} that \whp\ $\G(X,Y_j,p_1)$ contains a perfect matching
$M_j$. 
Now each perfect matching 
$M_j$ gives rise to an edge-coloured perfect matching $M^*_j$ of $G^*$
where $(x,y_i,x')$ gives rise to an edge $(x,x')$ of colour $y$. 
By symmetry, these matchings are uniformly random and they are independent by construction.
Also the edges have been randomly
coloured so that each colour appears exactly $r$ times. Indeed to achieve such a random colouring
we can take any partition of the edge set of $M_1^*\cup M_2^*\cup\cdots\cup M_{2r}^*$ 
into $2r$ sets $S_1,S_2,\ldots,S_{2r}$
of size $m$ and then colour the edges by using random bijections from $Y_j\to S_j$ for $j=1,2,\ldots,2r$.

We apply Theorems \ref{th3} and \ref{th4} to finish the
proof. Here the event $\cP_n$ of \eqref{cP} can be defined: 
\begin{multline*}
\cP_n=\\\ 
\{\text{Almost every equitable edge
colouring by $n$ colours produces a rainbow Hamilton cycle}\}.
\end{multline*}
By equitable, we mean that each colour is used $r$ times. This completes the proof of Theorem~\ref{th1}.

{\bf Acknowledgement} I am grateful to Michael Krivelevich and Oleg Pikhurko for their comments. 

\end{document}